\documentclass[10pt,leqno,twoside]{amsart}
\newtheorem{theorem}{Theorem}[section]
\newtheorem{lemma}[theorem]{Lemma}
\newtheorem{proposition}[theorem]{Proposition}
\newtheorem{corollary}[theorem]{Corollary}
\theoremstyle{definition}

\newtheorem{remark}[theorem]{Remark}
\newtheorem{remarks}[theorem]{Remarks}

\usepackage{latexsym}
\usepackage{amsmath}
\usepackage{amssymb}
\usepackage{mathrsfs}
\usepackage{graphicx}
\usepackage{ifthen}
\numberwithin{equation}{section}

\usepackage{latexsym}
\usepackage{amsmath}
\usepackage{amssymb}
\newcommand{\cA}{{\mathcal A}}

\newcommand{\cD}{{\mathcal D}}

\newcommand{\cE}{{\mathcal E}}
\newcommand{\cB}{{\mathcal L}}
\newcommand{\cL}{{\mathcal L}}
\newcommand{\cSM}{{\mathcal {SM} }}
\newcommand{\C}{{\mathbb C}}

\newcommand{\EE}{{\mathbb E}}
\newcommand{\N}{{\mathbb N}}
\newcommand{\R}{{\mathbb R}}
\newcommand{\RR}{{\mathbb R}}
\newcommand{\PP}{{\mathbb P}}

\newcommand{\dive}{\mbox{div }}
\newcommand{\curl}{\mbox{curl }}
\newcommand{\nn}{\nonumber}

\begin{document}

\title[Dynamics of the Ericksen-Leslie Equations I]
{Dynamics of the Ericksen-Leslie Equations with General Leslie Stress I: The Incompressible Isotropic Case}

\author{Matthias Hieber}
\address{Technische Universit\"at Darmstadt\\
        Fachbereich Mathematik\\
        Schlossgarten-Strasse 7\\
        D-64289 Darmstadt, Germany, and \\
607 Benedum Engineering Hall, University of Pittsburgh, Pittsburgh, PA 15261, USA}  
\email{hieber@mathematik.tu-darmstadt.de}

\author{Jan Pr\"uss}
\address{Martin-Luther-Universit\"at Halle-Witten\-berg\\
         Institut f\"ur Mathematik \\
         Theodor-Lieser-Strasse 5\\
         D-06120 Halle, Germany}
\email{jan.pruess@mathematik.uni-halle.de}

\subjclass[2000]{35Q35, 76A15, 76D03, 35K59}
\keywords{Ericksen-Leslie equations, nematic liquid crystals, general Leslie stress, global strong well-posedness, convergence to equilibria}

\begin{abstract}
The Ericksen-Leslie model for nematic liquid crystals in a bounded domain with general Leslie and isotropic Ericksen stress tensor is studied in the case  of a non-isothermal  and
incompressible fluid. This system is shown to be locally strongly well-posed in the $L_p$-setting, and a dynamical theory is developed.
The equilibria are identified and shown to be normally stable. In particular, a local solution extends to a unique, global strong solution provided the initial data are close to an
equilibrium or the solution is eventually bounded in the topology of the natural state manifold. In this case, the solution converges exponentially to an equilibrium, in the topology
of the state manifold. The above results are proven {\em without} any structural assumptions on the Leslie coefficients and in particular {\em without} assuming Parodi's relation.
\end{abstract}

\maketitle

\section{Introduction: The Ericksen-Leslie model with general Leslie stress}
In their pioneering articles, Ericksen \cite{Eri62} and Leslie \cite{Les68}  developed a  continuum theory for the flow of  nematic liquid crystals.  Their theory models nematic
liquid crystal flow from a hydrodynamical point of view and  describes the evolution of the underlying system under the influence of the velocity $u$ of the fluid and the
orientation configuration $d$ of rod-like liquid crystals. We already observe here that the modulus $|d|_2$ of the director field $d$ must equal 1 pointwise,
as $d$ represents a direction field. For more information see e.g. \cite{Cha92,DGP95,Ste04,ASS04} or \cite{LW14b}, \cite{HP15b}.
The original derivation \cite{Eri62,Les68} is based on the conservation  laws for mass, linear and angular momentums
as well as on constitutive relations given by Leslie in \cite{Les68}.


Following arguments from thermodynamics and employing the entropy principle, we proposed in \cite{HP15a}, \cite{HP15b} thermodynamically consistent models of Ericksen-Leslie type, 
even in the case of compressible fluids. Let us emphasize that these  models {\em contain the classical Ericksen-Leslie model in its general form as a special case}.

A related class of models also dealing with  the non-isothermal situation was presented by Feireisl, Rocca and Schimperna \cite{FRS11} as well as by Feireisl, Fr\'emond, Rocca and Schimperna
in \cite{FFRS12}. Their models  include stretching as well as rotational terms and are consistent with the fundamental laws of thermodynamics. The equation for the director $d$ is, however,
given in the penalized form, which does not seem to be physical. They show that the presence of the term $|\nabla d|_2^2$ in the internal energy as well as the stretching term
$d \cdot \nabla u$  give rise, in order to respect the laws of thermodynamics, to two new non dissipative contributions in the stress tensor $S$ and in the flux $q$. It is interesting
to note that these two new contributions coincide with the extra terms derived by Sun and Liu \cite{SL09} by different methods. It seems that these extra terms are non physical and arise there
for  purely mathematical reasons.


Given a bounded domain $\Omega \subset \R^n$, $n \geq 2$, with smooth enough boundary,  the general Ericksen-Leslie model in the  non-isothermal situation derived as in
\cite{HP15a} and \cite{HP15b} reads as
\begin{align}\label{eq:elgeneral}
\left\{
\begin{array}{rll}
\partial_t \rho +{\rm div}(\rho u)&=0\quad &\mbox{in } \Omega,\\
\rho\cD_tu +\nabla \pi &= {\rm div}\, S \quad &\mbox{in } \Omega,\\
\rho\cD_t\epsilon +{\rm div}\, q &= S:\nabla u -\pi{\rm div}\, u+{\rm div}(\rho\partial_{\nabla d}\psi\cD_td)\quad &\mbox{in } \Omega,\\
\gamma \cD_td-\mu_V Vd&= P_d\big( {\rm div}(\rho\frac{\partial\psi}{\partial\nabla d})-\rho\nabla_d\psi\big)+ \mu_D P_d Dd \quad &\mbox{in } \Omega,\\
u=0,\quad q\cdot \nu&=0 \quad  &\mbox{on } \partial\Omega,\\
\rho(0)=\rho_0,\quad u(0)=u_0,\quad \theta(0)&=\theta_0,\quad d(0)=d_0\quad &\mbox{in } \Omega.
\end{array}\right.
\end{align}
Here the unknown variables $\rho, u$ denote the density and velocity of the fluid, respectively, $\epsilon$ the
internal energy and $d$ the so called director, which - we recall - must have modulus 1.  Moreover, $D = \frac{1}{2}(\nabla u]^T+ \nabla u)$ denotes the rate of deformation tensor, the vorticity tensor $V$ defined by
$V=\frac{1}{2}(\nabla u -[\nabla u]^T)$ is skew-symmetric, and $q$ denotes the heat flux.  By  $\cD_t=\partial_t +u\cdot\nabla$ we denote the Lagrangian derivative and $P_d$ is defined as $P_d=I-d\otimes d$.

Note that the condition $|d|_2=1$ is preserved by smooth solutions, as $P_d d=0$ as well as $(Vd|d)=0$, hence $\cD_t|d|^2_2=0$.

In addition, $\psi$ denotes the  free energy, which, following Oseen  \cite{Ose33} and Frank \cite{Fra58}, see also Virga \cite{Vir94}, is given by the Oseen-Frank functional $\psi^F(d,\nabla d)$
as
\begin{equation}\label{eq:psif}
\psi^F(d,\nabla d) := k_1(\dive d)^2 + k_2 (d \cdot \curl d)^2 + k_3|d \times \curl d|^2 + (k_2+k_4)(tr(\nabla d)^2- (\dive d)^2),
\end{equation}
where $k_1,\ldots,k_4$ are the so called {\em Frank coefficients}, which may depend on $\rho$ and $\theta$.

These equations have to be supplemented by the thermodynamical laws
\begin{align}\label{eq:thermo}
\epsilon = \psi +\theta \eta, \quad \eta = -\partial_\theta \psi, \quad \kappa = \partial_\theta \epsilon,\quad \pi =\rho^2 \partial_\rho\psi,
\end{align}
and by the constitutive laws
\begin{align}\label{eq:const}
\left\{
\begin{array}{lll}
S &= S_N + S_E + S_L^{stretch} + S_L^{diss}, \\
S_N&= 2\mu_s D + \mu_b {\rm div}\, u \, I, \\
S_E&=  -\rho \frac{\partial \psi}{\partial\nabla d}[\nabla d]^{\sf T},\\
S_L^{stretch}& = \frac{\mu_D+\mu_V}{2\gamma}{\sf n} \otimes d + \frac{\mu_D-\mu_V}{2\gamma} d \otimes {\sf n},\quad
 {\sf n}= \mu_V Vd +\mu_D P_dDd-\gamma\cD_td,\\
S_L^{diss}& =\frac{\mu_P}{\gamma} ({\sf n}\otimes d + d\otimes {\sf n}) + \frac{\gamma\mu_L+\mu_P^2}{2\gamma}(P_dDd\otimes d + d\otimes P_dDd) +\mu_0 (Dd|d)d\otimes d,\\
\end{array}\right
.\end{align}
and
\begin{equation}\label{eq:q1}
q = -\tilde\alpha_0 \nabla\theta -\tilde\alpha_1(d|\nabla\theta)d.
\end{equation}
Here all coefficients $\mu_j,\tilde\alpha_j$ and $\gamma$ are functions of $\rho,\theta,d,\nabla d$.
For thermodynamical consistency we require the conditions
\begin{align}\label{eq:pos1}
\mu_s\geq0,\quad 2\mu_s+n\mu_b\geq0,\quad \tilde\alpha_0\geq0,\quad \tilde\alpha_0+\tilde\alpha_1\geq0,\quad \mu_0,\mu_L\geq0,\quad \gamma>0.
\end{align}
We also note that the natural boundary condition at $\partial\Omega$ for $d$ becomes
\begin{equation}
 \nu_i\nabla_{\partial_i d}\psi = 0\quad \mbox{ on }\partial\Omega.
\end{equation}
Observe that this condition is {\em fully nonlinear}, in general. Physically, it means that the boundary does not interact with the director field. 
Otherwise one would have  to model such interactions and it seems to be unclear  whether this could  be done in a physically consistent way by simply imposing Dirichlet boundary conditions. 
For this reason we employ the {\em Neumann condition} for $d$ throughout this paper. Actually, we can prove local well-posedness also in the case of Dirichlet boundary conditions, but the 
set of equilibria becomes more complicated in this case. Our results on stability and long time behaviour are only valid for constant equilibria.

In the case of {\em isotropic elasticy} with constant density and temperature one has $k_1=k_2=k_3=1$ and $k_4=0$ and so the Oseen-Frank energy reduces to the Dirichlet energy, i.e.
$$
\psi(d,\nabla d):= \psi^F(d,\nabla d) = \frac{1}{2}|\nabla d|^2,
$$
and  thus
$
\dive ( \frac{\partial \psi^F}{\partial(\partial \nabla d)}) = \Delta d$.
Then the Ericksen stress tensor simplifies to
\begin{equation}\label{eq:ericksenstress}
S_E=  -\lambda \nabla d[\nabla d]^{\sf T},
\end{equation}
where $\lambda = \rho\partial_\tau\psi$,
and the natural boundary condition at $\partial\Omega$ for $d$ becomes the Neumann condition $\partial_\nu d = 0$ on $\partial\Omega$.

It is the aim of this article to investigate the above Ericksen-Leslie system analytically, in the case of isotropic elasticity and for incompressible as well as  compressible fluids.
To this end, Part I of this article will concentrate on the case of incompressible fluids, whereas Part II will investigate the compressible case.
Here {\em incompressibility} means that the density $\rho$ is constant and {\em isotropy} means that the free energy $\psi$ is a function of $\varrho$,  $\theta$ and $\tau=|\nabla d|_2^2/2$, only.

For the convenience of the reader, we rewrite the model in the incompressible and isotropic case, which then reads as follows.
\begin{align}\label{eq:el}
\left\{
\begin{array}{rll}
\rho\cD_tu +\nabla \pi &= {\rm div}\, S \quad &\mbox{in } \Omega,\\
{\rm div}\, u&=0\quad &\mbox{in } \Omega,\\
\rho\cD_t\epsilon +{\rm div}\, q &= S:\nabla u+ {\rm div}(\lambda\nabla d\cD_t d)\quad &\mbox{in } \Omega,\\
\gamma\cD_td -\mu_V Vd- {\rm div}[ \lambda\nabla]d&=\lambda|\nabla d|^2d  + \mu_D P_d Dd\quad &\mbox{in } \Omega,\\
u=0,\quad q\cdot \nu=0,\quad \partial_\nu d&=0 \quad &\mbox{on } \partial\Omega,\\
\rho(0)=\rho_0,\quad u(0)=u_0,\quad \theta(0)&=\theta_0,\quad d(0)=d_0\quad &\mbox{in } \Omega.
\end{array}\right.
\end{align}
These equations have to be supplemented by the thermodynamical laws for the internal energy $\epsilon$, entropy $\eta$, heat capacity $\kappa$
given in \eqref{eq:thermo}, by Ericksen's tension $\lambda =\rho\partial_\tau\lambda$, and by the constitutive laws \eqref{eq:const} as above, with Ericksen stress of the form \eqref{eq:ericksenstress}, and with $q$ satifying
\begin{equation}\label{eq:q2}
q = -\alpha \nabla\theta.
\end{equation}
Note that all coefficients $\mu_j,\alpha$ and $\gamma$ are functions of $\theta$ and $\tau$, by the principle of equi-presence.

For further purposes, it  is convenient to write the equation for the internal energy as an equation for the temperature $\theta$. It reads as
$$
\rho\kappa\cD_t \theta +{\rm div}\, q = S:\nabla u + {\rm div}(\lambda \nabla)d\cdot \cD_t d + (\theta \partial_\theta\lambda) \nabla d \nabla\cD_t d.
$$
Observe the appearance of unusual third order terms due to the presence of $\cD_t d$ in the Leslie stress $S_L$ as well as in the last term of the energy balance.  This alludes a peculiarity of the system, which has to be overcome in the analysis.

Let us emphasize that in the case $\mu_V=\gamma$, our parameters $\mu_s,\mu_0,\mu_V,\mu_D,\mu_P,\mu_L$ are in one-to-one correspondence to the celebrated
Leslie parameters $\alpha_1,\ldots,\alpha_6$ given in the Leslie stress $\sigma_L$ defined by
\begin{equation}\label{eq:leslie}
\sigma_L  := \alpha_1 (d^TDd)d\otimes d +\alpha_2N \otimes d + \alpha_3 d\otimes N + \alpha_4 D + \alpha_5(Dd)\otimes d + \alpha_6 d \otimes (Dd),
\end{equation}
where $D$ denotes the deformation tensor as above and
$$
N:=N(u,d):=\partial_td + (u \cdot \nabla)d - Vd,
$$
with $V$ as above. This shows that our model \eqref{eq:el}, \eqref{eq:thermo}, \eqref{eq:const} {\em contains the classical isothermal and isotropic
Ericksen-Leslie model} given by
\begin{align}\label{eq:el-isotropic}
\left\{
\begin{array}{rlll}
 u_t+ (u \cdot \nabla) u -\Delta u  +\nabla \pi&\!=\!&   -\dive(\nabla d[\nabla d]^T) + \dive \sigma_L  & \text{in } (0,T) \times \Omega,  \\
 \dive u &\!=\!& 0 & \text{in } (0,T) \times \Omega, \\
 d_t + u \nabla d - Vd +  \frac{\lambda_2}{\lambda_1}Dd  &\!=\!& -\frac{1}{\lambda_1}(\Delta d + |\nabla d|^2d) + \frac{\lambda_2}{\lambda_1}(d^TDd)d & \text{in } (0,T) \times \Omega,\\
  (u,\partial_\nu d) &\!=\!& (0,0)& \text{on } (0,T) \times \partial \Omega, \\
  (u,d)_{\vert t=0} &\!=\!& (u_0,d_0) & \text{in }   \Omega. \end{array}\right.
\end{align}
{\em as a special case}; here we have $\lambda_1= -\gamma/\lambda$ and $\lambda_2 = \mu_D/\lambda$. Note that in the simplified model no third order terms appear, which considerably simplifies this problem.

It was shown by Parodi \cite{Par70} in 1970 that, {\em assuming Onsager's reciprocal relation}, one is lead to the relation
\begin{equation}\label{eq:parodi}
\alpha_2 + \alpha_3 = \alpha_6 - \alpha_5,
\end{equation}
where the coefficients $\alpha_j$ denote the Leslie coefficients introduced above.

The analysis of the Ericksen-Leslie system began by the pioneering work of Lin \cite{Lin89}, \cite{Lin91} and  Lin and Liu \cite{LL95}, \cite{LL00}, who introduced and studied
the nowadays called {\em isothermal simplified model}. They studied the situation where the nonlinearity in the equation for $d$ is replaced by a Ginzburg-Landau energy functional
and proved the existence of global weak solutions under suitable assumptions on the intial data. Wang proved in \cite{Wan11} global well-posedness for the simplified  system for
initial data being small in $BMO^{-1}\times BMO$ in the case of a whole space $\Omega=\R^n$ by combining techniques of Koch and Tataru with methods from
harmonic maps.

Concerning the situation of {\em bounded domains}, a rather complete understanding of the well-posedness as well as the  dynamics of the simplified system  subject to
Neumann conditions for $d$ was obtained in \cite{HNPS14}.  First results on the existence of global weak solutions to the  simplified system subject to Dirichlet boundary conditions 
in two dimensions go back to F.Lin, J. Lin and Wang \cite{LLW10}. Recently, considering the simplified system in three dimensions subject to Dirichlet conditions $d=d_b$ on $\partial\Omega$, 
Huang, Lin, Liu and Wang \cite{HLLW15} constructed  examples of small initial data  for which one has finite time blow up of $(u,d)$.

It seems that Coutard and Shkoller \cite{CS01} were the first who considered so called
{\em stretching terms} analytically in the equation for $d$. More precisely, they replaced the equation for $d$ in the {\em simplified model}
by a Ginzburg-Landau type approximation {\em including stretching} of the form
\begin{equation}\label{eq:introstretch}
 \gamma(\partial_td  + u \cdot \nabla d  + d \cdot \nabla u)= \Delta d - \frac{1}{\varepsilon^2}(|d|_2^2 -1)d \quad  \text{in } (0,T) \times \Omega.
\end{equation}
They proved local well-posedness for \eqref{eq:introstretch} as well as a  global existence result for small data.
Note, however, that in this case the presence of the stretching term $d \cdot \nabla u$ causes loss of total energy balance and, moreover, the
condition $|d|_2= 1$  in  $(0,T) \times \Omega$, is not preserved anymore.

For recent results on the general Ericksen-Leslie model with vanishing Leslie  and general Ericksen stress we refer to the articles
\cite{HLX13} and \cite{MLG14}. For results including stretching terms for $d$, we refer to the articles
\cite{LLW10}, \cite{HLW13}, \cite{WZZ13}, \cite{WXL13},  and \cite{HP15a,HP15b},  which contain well-posedness criteria  for the
general system under various assumptions on the Leslie coefficients.



The main new idea in investigating the general Ericksen-Leslie equations analytically in the strong sense is, similiarly to the situation of the simplified and isothermal
system treated in \cite{HNPS14}, to regard them  as a quasilinear parabolic evolution equation.
Restricting ourselves in Part I to the case of {\em incompressible fluids} and to the case of isotropic elasticity,
we present in the following a rather complete dynamic theory for the underlying equations. It seems to be the first well-posedness result for the {\em Ericksen-Leslie equations
dealing with general Leslie stress $S_L$ without assuming additional conditions on the Leslie coefficients}.

The results given in the three main theorems  below answer all questions concerning well-posedness, stability and longtime  behaviour for the {\em general} Ericksen-Leslie system subject 
to Neumann boundary conditions for $d$ in a very  satisfactory way.
It is proved by means of techniques involving maximal $L_p$-regularity and quasilinear parabolic evolution equations. For these methods, we refer to the booklet  by
Denk, Hieber, and Pr\"{u}ss \cite{DHP03}, the articles \cite{Pru03}, \cite{PSZ09}, \cite{KPW10}, \cite{LPW14}  and to the  monograph by Pr\"uss and Simonett \cite{PS16}. For the convenience of  the reader we have summarized the relevant results from these papers in Section3.

Let us also emphasize that for obtaining our well-posedeness results in the strong sense  {\em no structural conditions} on the Leslie coefficients are imposed and that in
particular {\em Parodi's relation}  \eqref{eq:parodi} on the Leslie coeffcients is {\em not} being assumed.

Moreover, the equilibria of the system have been identified in our recent paper [14] - which are zero velocities and constant density, temperature and director - and there it also has been proved that
these are thermodynamically stable. The negative total entropy has been shown to be a strict Lyapunov functional; in particular, the model is thermodynamically consistent.


We further mention at this point the series of articles \cite{FRS11},\cite{FFRS12},\cite{FSRZ12},\cite{FSRZ13}, in which it is shown that their particular
systems admit a global weak solution for a natural class of initial data.

For more information on modeling and analysis of the Ericksen-Leslie system, we refer e.g. to \cite{Cha92}, \cite{DGP95} and to the survey articles  \cite{HP15b} and \cite{LW14b}.

\section{Thermodynamical Stability and Consistency}
In this short section we recall from \cite{HP15a} and \cite{HP15b} that the above model \eqref{eq:el}, \eqref{eq:thermo}, \eqref{eq:const}, \eqref{eq:q2} has the following thermodynamical
properties. To this end, we introduce the following \\
{\bf Assumption (P)}: \\
\begin{align}\label{eq:pos2}
\mu_s>0,\quad \alpha>0,\quad \mu_0,\mu_L\geq0,\quad \kappa,\gamma>0,\quad \lambda, \lambda +2\tau \partial_\tau\lambda>0.
\end{align}

\begin{theorem}(\cite{HP15a}, Theorem 1). \label{prop:equi}
Assume that condition (P) holds. Then the  incompressible and isotropic model  \eqref{eq:el}, \eqref{eq:thermo}, \eqref{eq:const}, \eqref{eq:q2} has the following properties:
\begin{itemize}
\item[i)] Along smooth solutions the total energy
$$
{\sf E}:= \int_\Omega [\frac{\rho}{2}|u|_2^2 +\rho\epsilon]dx
$$
is preserved.
\item[ii)] Along smooth solutions the total entropy
$${\sf N}:= \int_\Omega \rho\eta dx$$
is non-decreasing.
\item[iii)] The negative total entropy $-{\sf N}$ is a strict Lyapunov functional.
\item[iv)] The condition $|d|_2=1$ is preserved along smooth solutions.
\item[v)] The equilibria of the system are given by the set of constants
\begin{equation}\label{eq:equilibria}
\cE=\{(0,\theta_*,d_*):\, \theta_*\in(0,\infty), \; d_*\in\RR^n,\; |d_*|_2=1\}.
\end{equation}
Here $\theta_*$ is uniquely determined by  the identity
$$ \epsilon(\theta_*,0) = {\sf E}_0/\rho|\Omega|.$$
\item[vi)] The equilibria are the critical points of the total entropy with prescribed energy.
\item[vii)] The second variation of ${\sf N}$ with prescribed energy at an  equilibrium is negative semi-definite.
\end{itemize}
\end{theorem}

\section{Background: Quasilinear Parabolic Evolution Equations}

\noindent
In this section we briefly recall some results on abstract quasilinear parabolic problems 
    \begin{equation}\label{Int1}
    \dot{v}+A(v)v=F(v),\quad t>0,\quad v(0)=v_0,
    \end{equation}
which are employed in the proofs of our main theorems. These results are due to Pr\"uss \cite{Pru03}, Pr\"uss and Simonett \cite{PS04}, K\"ohne, Pr\"uss and Wilke \cite{KPW10}, and
 Pr\"uss, Simonett and Zacher \cite{PSZ09}; a convenient reference for this theory is the monograph by Pr\"uss and Simonett \cite{PS16}, Chapter 5.
 
Assume that $(A,F):V_\mu\to \cB(X_1,X_0)\times X_0$ and $v_0\in V_\mu$. Here the spaces $X_1,X_0$ are Banach spaces such that $X_1\hookrightarrow X_0$ with dense embedding and 
$V_\mu$ is an open subset of the real interpolation space
    $$X_{\gamma,\mu}:=(X_0,X_1)_{\mu-1/p,p},\quad \mu\in(1/p,1].$$
We are mainly interested in solutions $v$ of \eqref{Int1} having \emph{maximal} $L_{p}$-\emph{regularity}, i.e.
    $$v\in H_{p}^1(J;X_0)\cap L_{p}(J;X_1)=:\EE_1(J), \mbox{ where } J = (0,T).$$
The trace space of this class of functions is given by $X_\gamma:=X_{\gamma,1}$.
However, to see and exploit the effect of parabolic regularization in the $L_p$-framework it is also useful to consider solutions in the class of weighted spaces
    $$v\in H_{p,\mu}^1(J;X_0)\cap L_{p,\mu}(J;X_1)=:\EE_{1,\mu}(J),\quad \mbox{which means } \; t^{1-\mu}v\in \EE_1(J). $$
The trace space for this class of weighted spaces is given by $X_{\gamma,\mu}$.
In our approach it is crucial to know that the operators $A(v)$ have the property of {\em maximal $L_{p}$-regularity}. Recall that an operator
$A_0$ in $X_0$ with domain $X_1$ has maximal $L_p$-regularity, if the linear problem
$$ \dot{v} +A_0v =f, \quad t\in J,\; v(0)=0,$$
admits a unique solution $v\in \EE_1(J)$, for any given $f\in L_p(J;X_0)=:\EE_0(J)$. It has been proved in \cite{PS04} that in this case maximal regularity also holds in the weighted case. 
\begin{proposition}\label{LWPthm}
Let $p\in (1,\infty)$, $v_0\in V_\mu$ be given and suppose that $(A,F)$ satisfies 
    \begin{equation}\label{LWP2}
    (A,F)\in C^{1}(V_\mu;\cB(X_1,X_0)\times X_0),
    \end{equation}
 for some $\mu\in (1/p,1]$. Assume in addition that $A(v_0)$ has maximal $L_p$-regularity.

Then there exist $a=a(v_0)>0$ and $r=r(v_0)>0$ with $\bar{B}_{X_{\gamma,\mu}}(v_0,r)\subset V_\mu$ such that problem \eqref{Int1} has a unique solution
    $$v=v(\cdot,v_1)\in \EE_{1,\mu}(0,a)\cap C([0,T];V_\mu),$$
on $[0,a]$, for any initial value $v_1\in \bar{B}_{X_{\gamma,\mu}}(v_0,r)$. In addition, 
$$ t\partial_t v \in \EE_{1,\mu}(0,a),$$
in particular, for each $\delta\in(0,a)$ we have
$$v\in H^2_p((\delta,a);X_0)\cap H^1_p((\delta,a);X_1)\hookrightarrow C^1([\delta,a];X_\gamma)\cap C^{1-1/p}([\delta,a];X_1),$$
i.e.\ the solution regularizes instantly.
\end{proposition}
\noindent
The next result provides information about the continuation of local solutions.

\begin{corollary}\label{LWPcor1}
Let the assumptions of Theorem \ref{LWPthm} be satisfied and assume that $A(v)$ has maximal $L_p$-regularity for all $v\in V_\mu$. Then the solution $v$ of \eqref{Int1} has a maximal 
interval of existence $J(v_0)=[0,t_+(v_0))$, which is characterized by the following alternatives:
\smallskip\\
{ (i)} Global existence: $t_+(v_0)=\infty$;  \smallskip\\
{(ii)}  $\liminf_{t\to t_+(v_0)} {\rm dist}_{X_{\gamma,\mu}}(v(t),\partial V_\mu)=0$; \smallskip\\
{(iii)} $\lim_{t\to t_+(v_0)} v(t)$ does not exist in $X_{\gamma,\mu}$.
\end{corollary}

\noindent
Next we assume that there is an open
set $V\subset X_\gamma$ such that
\begin{equation}
\label{AF}
(A,F)\in C^1(V,\cB(X_1,X_0)\times X_0).
\end{equation}
Let $ \cE\subset V\cap X_1$ denote the set of  equilibrium solutions of (\ref{Int1}), which means that
$$
v\in\cE \quad \mbox{ if and only if }\quad v\in V\cap X_1 \mbox{ and } A(v)v=F(v).
$$
Given an  element $v_*\in\cE$,  we assume that $v_*$ is
contained in an $m$-dimensional {\em manifold of equilibria}. This means that there
is an open subset $U\subset\R^m$, $0\in U$, and a $C^1$-function $\Psi:U\rightarrow X_1$,  such that
\begin{equation}
\label{manifold}
\begin{aligned}
& \bullet\
\text{$\Psi(U)\subset \cE$ and $\Psi(0)=v_*$,} \\
& \bullet\
 \text{the rank of $\Psi^\prime(0)$ equals $m$, and} \\
& \bullet\
\text{$A(\Psi(\zeta))\Psi(\zeta)=F(\Psi(\zeta)),\quad \zeta\in U.$}
\end{aligned}
\end{equation}
We suppose that
the operator $A(v_*)$ has the property of maximal $L_p$-regularity, and define the full linearization of \eqref{Int1} at $v_*$ by
\begin{equation}
\label{A0}
A_0w=A(v_*)w+(A^\prime(u_*)w)v_*-F^\prime(v_*)w\quad\hbox{for }
w\in X_1.
\end{equation}
\noindent
After these preparations we can state the following result on convergence of  solutions starting near $v_*$ which is called the {\em generalized principle of linearized stability}.
\begin{proposition}\label{th:1} 
Let $1<p<\infty$. Suppose $v_*\in V\cap X_1$ is an
equilibrium of (\ref{Int1}), and suppose that the functions
$(A,F)$ satisfy \eqref{AF}. Suppose further that $A(v_*)$ has the
property of maximal $L_p$-regularity and let $A_0$ be defined in
\eqref{A0}.
Suppose that $v_*$ is {\rm normally stable},
which means 
\begin{description}
\item[{\rm i)}] near $v_*$ the set of equilibria $\cE$ is a $C^1$-manifold in $X_1$ of dimension $m\in\N$,
\vspace{-1mm}
\item[{\rm ii)}] \, the tangent space for $\cE$ at $v_*$ is isomorphic to ${\sf N}(A_0)$,
\vspace{-1mm}
\item[{\rm iii)}] \, $0$ is a {\em semi-simple eigenvalue} of
$A_0$, i.e.\ $ {\sf N}(A_0)\oplus {\sf R}(A_0)=X_0$,
\vspace{-1mm}
\item[{\rm iv)}] \, $\sigma(A_0)\setminus\{0\}\subset \C_+=\{z\in\C:\, {\rm Re}\, z>0\}$.
\end{description}
Then $v_*$ is stable in $X_\gamma$, and there exists $\delta>0$ such
that the unique solution $v$ of (\ref{Int1}) with initial
value $v_0\in X_\gamma$ satisfying $|v_0-v_*|_{\gamma}\leq\delta$
exists on $\R_+$ and converges at an exponential rate in $X_\gamma$
to some $v_\infty\in\cE$ as $t\rightarrow\infty$.
\end{proposition}
\noindent
The next result contains information on bounded solutions in the presence of compact embeddings and of a strict Lyapunov functional.

\begin{proposition}\label{LTBthm}
Let $p\in (1,\infty)$, $\mu\in (1/p,1)$, $\bar{\mu}\in (\mu,1]$, with $V_\mu\subset X_{\gamma,\mu}$ open.
Assume that $(A,F)\in C^1(V_\mu;\cB(X_1,X_0)\times X_0)$, and that the embedding
$X_{\gamma,\bar{\mu}}\hookrightarrow X_{\gamma,\mu}$ is compact. Suppose furthermore that $v$ is a maximal solution which is bounded in $X_{\gamma,\bar{\mu}}$ and satisfies
\begin{equation}\label{GDC}
{\rm dist}_{X_{\gamma,\mu}}(v(t),\partial V_{\mu})\ge\eta>0,\; \mbox{ for all }t\geq 0.
\end{equation}
Suppose that $\Phi\in C(V_{\mu}\cap X_\gamma;\R)$ is a strict Lyapunov functional for \eqref{Int1}, which means that $\Phi$ is strictly decreasing along non-constant solutions.

Then $t_+(v_0)=\infty$, i.e.\ $v$ is a global solution of \eqref{Int1}. Its $\omega$-limit set $\omega_+(v_0)\subset \cE$  in $X_\gamma$ is nonempty, compact and connected.
If, in addition, there exists $v_*\in\omega_+(v_0)$ which is normally stable, then $\lim_{t\to\infty}v(t)=v_*$ in $X_\gamma$.
\end{proposition}

\section{Main Results}
In order to  formulate the main well-posedness result for the Ericksen-Leslie system in the incompressible and isotropic case and  to have access to the tools presented  in Section 3,
we introduce a functional analytic setting as follows. Denote the principal variable by  $v=(u,\theta,d)$ and
let us rewrite system  \eqref{eq:el}, \eqref{eq:thermo}, \eqref{eq:const}, \eqref{eq:q2} as a quasi-linear evolution equation of the form
\begin{equation} \label{eq:QLP}
\dot{v} +A(v) v = F(v),\quad t>0,\; v(0)=v_0,
\end{equation}
replacing $\cD_td$ appearing in the the equations for $u$ and $\theta$ by the equation for $d$. We also  apply the {\em Helmholtz projection} $\PP$ to the equation for $u$.
Then  $v$ belongs to the base space $X_0$ defined by
$$
X_0 :=L_{q,\sigma}(\Omega)\times L_q(\Omega;\R)\times H^1_q(\Omega;\R^n),
$$
where $1<p,q<\infty$ and $\sigma$ indicates solenoidal vector fields.
The regularity space will be
$$
X_1 :=\{ u\in H^2_q(\Omega;\RR^n)\cap L_{q,\sigma}(\Omega):\, u=0 \mbox{ {on} } \partial\Omega\}\times Y_1,
$$
with
$$
Y_1 := \{ (\theta,d)\in H^2_q(\Omega)\times H^3_q(\Omega;\RR^n):\, \partial_\nu \theta=\partial_\nu d =0 \mbox{ on } \partial\Omega\}.
$$
We consider solutions $v$ within the class
$$
v\in H^1_{p,\mu}(J;X_0)\cap L_{p,\mu}(J;X_1)=\EE_{1,\mu}(J),
$$
where $J=(0,a)$ with $0<a\leq\infty$ is an interval and $\mu\in (1/p,1]$.
The time trace space of this class is given by
\begin{equation}\label{eq:trace}
X_{\gamma,\mu} := \{ u\in B_{qp}^{2(\mu-1/p)}(\Omega)^n\cap L_{q,\sigma}(\Omega):\, u=0 \mbox{ on } \partial\Omega\}\times Y_{\gamma,\mu},
\end{equation}
where
$$
Y_{\gamma,\mu} = \{ (\theta,d)\in B^{2(\mu-1/p)}_{qp}(\Omega)\times B^{1+2(\mu-1/p)}_{qp}(\Omega;\R^n):\, \partial_\nu\theta=\partial_\nu d=0 \mbox{ on } \partial\Omega\},$$
whenever the boundary traces exist.
Note that
$$
X_{\gamma,\mu}\hookrightarrow B_{qp}^{2(\mu-1/p)}(\Omega)^{n+1}\times B_{qp}^{1+2(\mu-1/p)}(\Omega)^n\hookrightarrow C(\overline{\Omega})^{n+1}\times C^1(\overline{\Omega})^n,
$$
provided
\begin{equation}\label{eq:pq}
\frac{1}{p} +\frac{n}{2q}<\mu\leq 1.
\end{equation}
For brevity we set $X_\gamma:=X_{\gamma,1}$. Finally, the state manifold of the problem is defined by
$$
\cSM = \{ v\in X_{\gamma}:\, \theta(x)>0,\,  |d(x)|_2 =1 \mbox{ in } \Omega\}.
$$
We assume the following regularity assumptions on the parameter functions:

\vspace{.2cm}\noindent
{\bf Regularity assumption (R):} \\
The parameter functions are assumed to satisfy
\begin{equation}\label{eq:reg}
\mu_j, \alpha, \gamma\in C^2((0,\infty)\times[0,\infty)) \mbox{ for } j=S,V,D,P,L,0, \mbox{ and } \psi\in C^4((0,\infty)\times [0,\infty)).
\end{equation}

\vspace{.2cm}
The  fundamental well-posedness results regarding the general isotropic incompressible Ericksen-Leslie system reads  as follows.
\medskip

\begin{theorem}{\rm (Local Well-Posedness).} \label{thm:1} \\
Let $\Omega \subset \R^n$ be a bounded domain with boundary of class $C^{3-}$ and assume conditions (R) and (P).
Assume that  $J=(0,a)$, $1<p,q,<\infty$ and $\mu \in (1/p,1]$ are subject to \eqref{eq:pq} and  $v_0 \in X_{\gamma,\mu}$. Then for some $a=a(v_0)>0$, there is a unique solution
$$
v\in H^1_{p,\mu}(J,X_0)\cap L_{p,\mu}(J;X_1),
$$
of \eqref{eq:QLP}, i.e. of \eqref{eq:el}, \eqref{eq:thermo}, \eqref{eq:const}, \eqref{eq:q2}  on $J$. Moreover,
$$
v\in C([0,a];X_{\gamma,\mu})\cap C((0,a];X_\gamma),
$$
i.e.\ the solution regularizes instantly in time. It depends continuously on $v_0$ and exists on a maximal time interval $J(v_0) = [0,t^+(v_0))$.
Moreover,
\begin{align*}
t\partial_t v \in  H^1_{p,\mu}(J;X_0)\cap L_{p,\mu}(J;X_1),\quad a<t^+(v_0).
\end{align*}
Furthermore,
$$
|d(\cdot,\cdot)|_2 \equiv 1, \quad {\sf E}(t)\equiv {\sf E}_0, \quad  t \in J,
$$
and $-N$ is a strict Lyapunov functional.
In addition,  problem  \eqref{eq:QLP}   generates a local semi-flow in its natural state manifold $\cSM$.
\end{theorem}

Below, we denote by $\bar{\cE}$ the set
$$\bar{\cE} =\{ (0,\theta_*,d_*): \theta_*>0,\, d_*\in\RR^n\}$$
of constant equilibria of the system when ignoring the constraints $|d|_2=1$ and ${\sf E}={\sf E}_0$.  The next result concerns stability of equilibria.

\begin{theorem}{\rm (Stability of Equilibria).} \label{thm:2}\\
Assume conditions (R) and (P). Then any equilibrium $v_*\in \overline{\cE}$ of \eqref{eq:QLP} is stable in $X_\gamma$. Moreover, for each $v_*\in \overline{\cE}$
there is $\delta > 0$ such that if $|v_0 -v_*|_{X_{\gamma,\mu}} \leq \delta$, then the solution $v$ of \eqref{eq:QLP}  with initial value $v_0$ exists
globally in time and converges at an
exponential rate in $X_\gamma$ to some $v_\infty\in \overline{\cE}$.
\end{theorem}
\noindent
The third result concerns global existence and convergence of solutions to equilibria in the topology of the state manifold $\cSM$.

\begin{theorem}{\rm (Long-Time Behaviour).}\label{thm:3}\\
Assume conditions (R) and (P) and let $v$ be the solution of equation \eqref{eq:QLP} with $v_0 \in \cSM$. Then the following assertions hold. \\
{\bf a)} Suppose that for some $\bar{\mu}\in ( 1/p+n/2q+1/2,1]$ we have
\begin{equation}\label{eq:longtime}
\sup_{t \in (0,t^+(v_0))} [|v(t)|_{X_{\gamma,\mu}} + |1/\theta(t)|_{L_\infty}]<\infty.
\end{equation}
Then $t^+(v_0)=\infty$ and $v$ is a global solution.\\
{\bf b)} If $v$ is a global solution, bounded in $X_{\gamma,\bar{\mu}}$ and with $1/\theta$ bounded, then $v$ converges exponentially in $\cSM$ to an
equilibrium $v_\infty\in \mathcal E$ of \eqref{eq:QLP}   as $t\to \infty$.\\
{\bf c)} If $v$ is global solution of \eqref{eq:QLP} which converges to an equilibrium in $\cSM$, then \eqref{eq:longtime} valid.
\end{theorem}

\begin{remark}
Let us emphasize that the above theorem holds true {\em without any structural assumptions} besides condition {(P)} on the Leslie coefficients. In particular, the above
well-posedness results hold true {\em without} assuming Parodi's relation \eqref{eq:parodi}.
\end{remark}

\begin{remarks}
{a)} Wu, Xu and Liu considered in \cite{WXL13} the {\em isothermal penalized} Ericksen-Leslie model and gave a  formal physical derivation of the Ericksen-Leslie model based on an 
energy variational approach assuming Parodi's relation. Then they prove that, under certain assumptions on the data and the Leslie coefficients, the isothermal penalized Ericksen-Leslie 
system admits a unique, global solution provided the viscosity is large enough and study as well it longtime behaviour.  Moreover, {\em assuming Parodi's relation}, but not largeness of the 
viscosity, they show global well-posedness and Lypunov stability for the {\em penalized} Ericksen-Leslie system near local energy minimizers. \\
{b)} Wang, P. Zhang and Z. Zhang \cite{WZZ13} proved local well-posedness of the {\em  isothermal} general Ericksen-Leslie system as well as global well-posedness for small initial data
under various conditions on the Leslie coefficients, which ensure that the energy of the system is dissipated.\\
{c)} It is interesting to compare our above results with a recent result due to Huang, Lin, Liu and Wang \cite{HLLW15}, where they considered the simplified system subject to  
Dirichlet boundary conditions $d=d_b$ on $\partial\Omega$ and where they constructed examples of small initial data  for which one has finite time blow up of $(u,d)$ for the solution 
of the  simplified system. 
\end{remarks}

\section{Maximal $L^p$-Regularity of the Linearization}
The main task to apply the results in Section 3 is to establish maximal $L_p$-regularity of the linearized problem. To prove this,
we linearize equation \eqref{eq:el} at an initial value $v_0=(u_0,\theta_0,d_0)$ and drop all terms of lower order. This yields the principal linearization
\begin{align}\label{eq:prlin}
\left\{
\begin{array}{rll}
\cL_\pi(\partial_t,\nabla)v_\pi         &= f & \mbox{ in } J\times \Omega,\\
u=\partial_\nu\theta=\partial_\nu d &= 0 & \mbox{ on } J \times \partial\Omega,\\
u=\theta = d                       &= 0 & \mbox{ on } \{0\}\times\Omega.
\end{array}\right.
\end{align}
Here $J=(0,a)$, $v_\pi=(u,\pi,\theta,d)$ is the unknown, and $f=(f_u,f_\pi,f_\theta,f_d)$ are the given data. Denote the spatial co-variable by $\xi$ and by
$z$ that in time. Assume that $z \in \Sigma_\phi:=\{z \in \C\backslash\{0\}: 0 \leq |arg z| < \pi\}$. Then
the differential operator $L=\cL_\pi(\partial_t,\nabla)$ is defined via its symbol $\cL_\pi(z,i\xi)$, which is
\begin{equation}\label{princ-symb-pi}
\cL_\pi(z,i\xi) :=
\left[ \begin{array}{cccc}
M_u(z,\xi) & i\xi & 0& izR_1(\xi)^{\sf T}\\
i\xi^{\sf T}&0&0&0\\
0&0&m_\theta(z,\xi)& -iz \theta_0b a(\xi)\\
-iR_0(\xi)&0&-iba(\xi)& M_d(z,\xi)
\end{array}\right], \quad \xi \in \R^n, z \in \Sigma_\phi,
\end{equation}
with $b =\partial_\theta \lambda$, and $\lambda_1=\partial_\tau\lambda$. We also introduce the parabolic part of this symbol by dropping pressure gradient and divergence, i.e.
\begin{equation}\label{princ-symb}
\cL(z,i\xi) = \left[ \begin{array}{ccc}
M_u(z,\xi)  & 0& izR_1(\xi)^{\sf T}\\
0&m_\theta(z,\xi)& iz \theta_0b a(\xi)\\
-iR_0(\xi)&ib a(\xi)& M_d(z,\xi)
\end{array}\right], \quad \xi \in \R^n, z \in \Sigma_\phi.
\end{equation}
The entries of these matrices are given by
\begin{align*}
m_\theta(z,\xi) &:= \rho\kappa z +\alpha|\xi|^2,\\
a(\xi) &:=\xi\cdot\nabla d_0, \\
M_d(z,\xi) &:= \gamma z + \lambda|\xi|^2 +\lambda_1 a(\xi)\otimes a(\xi)= m_d(z,\xi) +\lambda_1 a(\xi)\otimes a(\xi),\\
R_0(\xi)&:= \frac{\mu_D +\mu_V}{2} P_0\xi\otimes d_0 + \frac{\mu_D -\mu_V}{2}(\xi| d_0) P_0,\\
R_1(\xi) &:=  (\frac{\mu_D +\mu_V}{2} +\mu_P)P_0\xi\otimes d_0 + (\frac{\mu_D -\mu_V}{2}+\mu_p)(\xi| d_0) P_0,\\
M_u(z,\xi) &:= \rho z + \mu_s|\xi|^2 + \mu_0(\xi|d_0)^2d_0\otimes d_0 +a_1(\xi|d_0)P_0\xi\otimes d_0 \\
& \quad + a_2 (\xi|d_0)^2 P_0 + a_3 |P_0\xi|^2 d_0\otimes d_0 + a_4 (\xi|d_0) d_0\otimes P_0\xi.
\end{align*}
Here $P_0=P_{d_0}=I-d_0\otimes d_0$, and $a_j$ are certain coefficients. Note that the above coefficients depend on $x$ through the dependence of the parameter functions on the
initial value $v_0(x)$. The maximal regularity result  for \eqref{eq:prlin} employed below  reads as follows.

\begin{theorem}\label{thm:mr}
Let $J=(0,a)$, $1<p,q<\infty$, and assume condition (R) and (P). Then  equation \eqref{eq:prlin} admits a unique solution $v_\pi=(u,\pi,\theta,d)$ satisfying
\begin{align*}
(u,\theta)&\in {_0H}^1_p(J;L_q(\Omega))^{n+1} \cap L_p(J; H^2_q(\Omega))^{n+1},\\
\pi&\in L_p(J; \dot{H}^1_q(\Omega)),\\
 d&\in {_0H}^1_p(J;H^1_q(\Omega))^{n} \cap L_p(J; H^3_q(\Omega))^{n},
\end{align*}
if and only if
\begin{align*}
\begin{array}{rl}
(f_u,f_\theta) &\in L_p(J;L_q(\Omega))^{n+1},\\
 f_d           &\in L_p(J;H^1_q(\Omega))^n, \\
f_\pi          &\in {_0H}^1_p(J;H^{-1}_q(\Omega))\cap L_p(J;H^1_q(\Omega)).\\
\end{array}
\end{align*}
Further, the solution map $f\mapsto v_\pi$ is continuous between the corresponding spaces.
\end{theorem}

Let us remark that if we replace $\partial_t $ by $\partial_t+\omega$, where $\omega>0$ is a sufficiently large constant, then the assertion of Theorem \ref{thm:mr}  holds true also for
$J=(0,\infty)$.

\vspace{.2cm}\noindent
{\em Proof}. We subdivide the proof into 5 steps.\\
{\em Step 1: The Principal Symbol with Constant Coefficients in $\Omega=\RR^n$.}\\
To extract the structure of $\cL$, we introduce the symbols
$$
R(\xi):= (\xi|d_0) P_0 + P_0\xi\otimes d_0,\quad R_\mu(\xi):= \mu_-(\xi|d_0) P_0 + \mu_+P_0\xi\otimes d_0,\quad \mu_\pm := \mu_D \pm \mu_V+\mu_P.
$$
Then $M_u$ simplifies to
$$
M_u= m_u + \mu_0(\xi|d_0)^2d_0\otimes d_0 + \frac{\mu_L}{4} R^{\sf T}R + \frac{1}{4\gamma} R_\mu^{\sf T}R_\mu + \frac{\mu_P\mu_V}{2\gamma}(\xi|d_0)(R-R^{\sf T}),
$$
and we also have
$$
R_1 = R_\mu -R_0 \mbox{ and } m_u(z,\xi) =\rho z +\mu_s|\xi|^2.
$$
Next, we set $v=(u,w)$, $v_\pi=(v,\pi,w)$ and $w=(\theta,d)$. Then, setting
$$
J={\rm diag}(I,1/\theta_0,zI),
$$
and  the second line of $\cL$ by $1/\theta_0$ as well as the last line with $\bar{z}$, we obtain the estimate
\begin{align*}
{\rm Re} (\cL v| Jv)& = {\rm Re}\, m_u |u|_2^2 + {\rm Re}\, m_\theta|\theta|^2 +  {\rm Re}\, z(\lambda_0|\xi|^2|d|_2^2 + \lambda_1|(a(\xi)|d)|^2)\\
& \quad + \frac{\mu_L}{4} |Ru|_2^2 + \frac{1}{4\gamma} |R_\mu u|_2^2 + {\rm Re}[ i z (d|R_\mu u) + \gamma|z|^2|d|_2^2\\
&\geq c [ {\rm Re}\, z( |u|_2^2 +|\theta|^2 + |\xi|^2|d|_2^2) + |\xi|^2(|u|_2^2+|\theta|^2) + (2\gamma|z||d|_2-|R_\mu u|_2)^2],
\end{align*}
provided
$$
\rho,\mu_s,\kappa,\gamma, \alpha, \lambda,\lambda +2\tau\partial_\tau\lambda>0 \mbox{ and } \mu_0,\mu_L\geq0.
$$
One could even relax the assumptions on $\mu_0$ and $\mu_L$ to $2\mu_s+\mu_0>0$ and $2\mu_s +\mu_L>0$, but we will not do this here.
This means that the symbol $\bar{J}\cL$ is accretive for ${\rm Re}\,z>0$, i.e. it is  strongly elliptic.

Let us emphasize  that we do not need any structural conditions on the coefficients $\mu_D,\mu_V,\mu_P, \partial_\theta\lambda$.

\vspace{.2cm}\noindent
{\em Step 2: Schur Reductions}. \\
In this step we perform  a Schur reduction  to reduce the above symbol to a symbol only for $u$.  To this end, we consider the subsystem for $w$, i.e.\  the equation
\begin{align*}
\left[\begin{array}{cc}
m_\theta(z,\xi) & -iz \theta_0b a(\xi)^{\sf T}\\
-ib a(\xi)& m_d(z,\xi) + \lambda_1 a(\xi)\otimes a(\xi)
\end{array}\right]
\left[\begin{array}{c}
\theta\\ d
\end{array}\right]= \left[ \begin{array}{c}
f_\theta\\ f_d + i R_0(\xi)u
\end{array}\right].
\end{align*}
To solve this system, we follow the strategy developed in \cite{HP15a} and introduce the new variable $\delta = (a(\xi)|d)$. Then, multiplying the second equation with
$a(\xi)$ we obtain the system
\begin{align*}
\left[\begin{array}{cc}
m_\theta(z,\xi) & -iz \theta_0b \\
-ib |a(\xi)|^2& m_d(z,\xi) + \lambda_1| a(\xi)|^2
\end{array}\right]
\left[\begin{array}{c}
\theta\\ \delta
\end{array}\right]= \left[ \begin{array}{c}
f_\theta\\ (f_d|a(\xi)) + i (R_0(\xi)u|a(\xi))
\end{array}\right].
\end{align*}
This system is easily solved to the result
\begin{align*}
\left[ \begin{array}{c} \theta\\
\delta\end{array}\right] =
\frac{1}{det(z,\xi)} \left[\begin{array}{cc}
m_d(z,\xi) + \lambda_1| a(\xi)|^2& iz \theta_0b \\
ib |a(\xi)|^2& m_\theta(z,\xi)
\end{array}\right]
\left[ \begin{array}{c}
f_\theta\\ (f_d + i R_0(\xi)u|a(\xi))
\end{array}\right],
\end{align*}
where
$$
det(z,\xi) = m_\theta(z,\xi)(m_d(z,\xi)+\lambda_1 |a(\xi)|^2) + z\theta_0b^2|a(\xi)|^2.
$$
Note that this symbol behaves like $(z+|\xi|^2)^2$ as soon as $\rho,\kappa,\lambda, \lambda+2\tau\partial_\tau\lambda>0$.

Knowing $\delta=(a(\xi)|d)$ and $\theta$, we are now able to determine $d$. As a result we obtain
$$
d = m_d^{-1} [ f_d +iR_0u + i b a(\xi)\theta - \lambda_1 a(\xi)\delta].
$$
Following the arguments given in \cite{HP15a}, we see that $\theta$ and $d$ belong to the right regularity classes,
whenever $f_\theta,f_d$ and $u$ are so.

In order to  extract the Schur complement for $u$, we set $f_\theta=f_d=0$ and compute $d$. This yields
\begin{equation}\label{eq:schurc}
d= i[\frac{1}{m_d}(I-a_0\otimes a_0)+ \frac{m_\theta}{det}a_0\otimes a_0]R_0 u=i[\frac{1}{m_d}P_{a_0}+ \frac{m_\theta}{det}Q_{a_0}]R_0 u.
\end{equation}
with $a_0(\xi)=a(\xi)/|a(\xi)|$ if $a(\xi)\neq0$ and $a_0(\xi)=0$ otherwise. This is the representation of $d$ needed for the Schur complement of $u$.

\vspace{.2cm}\noindent
{\em Step 3: The Generalized Stokes Symbol}.\\
We insert \eqref{eq:schurc} into the equation for $u$ to obtain the generalized Stokes symbol for $(u,\pi)$ and obtain
\begin{align}\label{eq:genstokes}
M(z,i\xi) &= M_u(z,\xi) - zR_1^{\sf T}(\xi)[\frac{1}{m_d(z,\xi)}P_{a_0}(\xi)+ \frac{m_\theta(z,\xi)}{det(z,\xi)}Q_a(\xi)]R_0(\xi)
\end{align}
As the Schur reduction preserves accretivity, even with the same accretivity constant, we obtain
$$
{\rm Re}(M(z,i\xi)u|u)\geq {\rm Re}\, m_u(z,\xi)|u|^2 =\big(\rho {\rm Re}\, z +\mu_s |\xi|^2\big)|u|^2.
$$
This  shows that $M$ is strongly elliptic. For this reason we may now apply the method developed by  Bothe and Pr\"uss \cite{BP07} or Pr{\"u}ss and Simonett \cite{PS16}, Section 7.1,
to prove maximal $L_p$-regularity of the resulting generalized Stokes problem. In these references we need to replace $\lambda$ by $z$ and the symbol $z +\cA(\xi)$ by $M(z,i\xi)$.
We will not do this here in detail and refer the reader to Section 7.1 of \cite{PS16} for this analysis.

\vspace{.2cm}\noindent
{\em Step 4: The Lopatinskii-Shapiro Condition}.\\
In order to guarantee  the solvability of the above problem in a half-space, we need  to replace the co-variable $\xi$ by the one-dimensional differential operator
$\xi-i\nu\partial_y$, where $(\xi|\nu)=0$. The Lopatinskii-Shapiro condition then means that the problem
\begin{align}\label{eq:lopshi}
\cL(z,i\xi+\nu\partial_y) v&=0,\quad y>0,\\
u(0)=\partial_y \theta=\partial_y d&=0,\nn
\end{align}
admits  only the zero solution in $L_2(\RR_+)^{2n+1}$, for all $(z,\xi)\neq(0,0)$.

In order to prove this condition, suppose that $v(y)$ is a solution of the ODE system \eqref{eq:lopshi}, which belongs to $L_2(\R_+)^{2n+1}$. Taking the inner product with
$v$ in $L_2(\R_+)$ , taking real parts, integrating by parts with respect to \ $y$ and employing the boundary conditions, we obtain the estimate
\begin{align*}
c{\rm Re}(\cL(z,i\xi+\nu\partial_y) v|v)_{L_2}&\geq {\rm Re}\, z [ |u|_{L_2}^2 +|\theta|_{L_2}^2 +|\xi|^2|d|_{L_2}^2+|\partial_y d|^2_{L_2} ]+|z|^2|d|_{L_2}^2\\
 & \quad + |\xi|^2(|u|_{L_2}^2 +|\theta|_{L_2}^2) + |\partial_yu|_{L_2}^2 + |\partial_y\theta|_{L_2}^2.
\end{align*}
This shows that the Lopatinskii-Shapiro is valid. Hence, to prove maximal $L_p$-regularity in the half space case, we may proceed in the following way. First we perform the
same Schur reductions as in Step 2 and as in \cite{HP15a}. This yields the unique existence of  $\theta$ and $d$ in the right regularity class. We then employ the half-space theory for
the generalized Stokes symbol $M$ by the methods in Bothe and Pr\"uss \cite{BP07} or Pr\"uss and Simonett \cite{PS16}, Section 7.2, to obtain maximal $L_p$-regularity for the
half-space case.

\vspace{.2cm}\noindent
{\em Step 5: General Domains and Variable Coefficients}.\\
The results of Step 3 and Step 4 extend by a perturbation argument to a bent half-space, and to the case of variable coefficients with small deviation from constant ones.
We then may apply a localization procedure to cover the case of general domains with smooth boundaries and variable coefficients. For details we refer at this point e.g.\ to
Sections 6.3 and 7.3. of the monograph \cite{PS16} by Pr\"uss and Simonett. This  completes the proof of Theorem \ref{thm:mr}.

\rightline{$\Box$}

\section{Proofs of the Main Results}
In this section we present the proofs of the above three main results. They are based on the theory of quasilinear parabolic evolution equations, see Section 3.

\vspace{.1cm}\noindent
{\em Proof of Theorem \ref{thm:1}}. \\
As already discussed above, we rewrite the system \eqref{eq:el}, \eqref{eq:thermo}, \eqref{eq:const}, \eqref{eq:q2} as a quasi-linear evolution equation of the form
\begin{equation} \label{NLCF}
\dot{v} +A(v) v = F(v),\quad t>0,\; v(0)=v_0,
\end{equation}
replacing $\cD_td$ appearing in the the equations for $u$ and $\theta$ by the equation for $d$.
Here $v=(u,\theta,d)$. We further apply the Helmholtz projection $\PP$ to the equation for $u$ and recall the  base space
$$
X_0= L_{q,\sigma}(\Omega)\times Y_0,
$$
with $Y_0=L_q(\Omega)\times H^1_q(\Omega;\R^n)$ as well as the regularity space $X_1$ as above, i.e.
$$
X_1 =\{ u\in H^2_q(\Omega;\RR^n)\cap L_{q,\sigma}(\Omega):\, u=0 \mbox{ {on} } \partial\Omega\}\times Y_1,
$$
with
$$
Y_1 = \{ (\theta,d)\in H^2_q(\Omega)\times H^3_q(\Omega;\RR^n):\, \partial_\nu \theta=\partial_\nu d =0 \mbox{ on } \partial\Omega\}.
$$
In order to prove local well-posedness of the system \eqref{eq:el}, \eqref{eq:thermo}, \eqref{eq:const}, \eqref{eq:q2} we may now resort to the abstract theory presented in Section 3. 

Note first that by  Theorem \ref{thm:mr} the quasi-linear part $A(v)$ has maximal $L_p$-regularity. A result by Pr{\"u}ss and Simonett \cite{PS04}, Theorem 2.4,  implies that $A(v)$
also admits maximal regularity in $L_{p,\mu}(J;X_0)$, hence also in the situation of time weights.
Recalling  the solution space
$\EE_\mu(J)= H^1_{p,\mu}(J;X_0)\cap L_{p,\mu}(J;X_1)$, we see that the time-trace space $X_{\gamma,\mu}$ of $\EE_\mu(J)$ is given as in \eqref{eq:trace}, and the embedding
\begin{equation}\label{eq:embedding}
X_{\gamma,\mu}\hookrightarrow B_{qp}^{2(\mu-1/p)}(\Omega)^{n+1}\times B_{qp}^{1+2(\mu-1/p)}(\Omega)^n\hookrightarrow C^1(\overline{\Omega})^{n+1}\times C^2(\overline{\Omega})^n,
\end{equation}
holds, provided
\begin{equation*}
\frac{1}{p} +\frac{n}{2q} +1/2<\mu\leq 1.
\end{equation*}
Here $B^s_{pq}(\Omega)$ denote as usual the Besov spaces; see e.g.\ Triebel \cite{Tri92}. Therefore, the mappings  $A$ and $F$ satisfy the assumptions 
of the local existence theorem Theorem 3.1, as well as of Corollary 3.2, hence we obtain local well-posedness for \eqref{NLCF} and strong solutions on a maximal time interval.

Even more, if only  $$\frac{1}{p} +\frac{n}{2q}<\mu\leq 1$$ holds, also the results in
LeCrone, Pr\"{u}ss and Wilke \cite{LPW14}, Theorem 2.1, apply, and we obtain local strong solutions if the initial values only satisfy
$$ (u_0,\theta_0,d_0)\in B^{2(\mu-1/p)}_{qp}(\Omega)^{n+1}\times B_{qp}^{1+2(\mu-1/p)}(\Omega)^n$$
plus compatibility conditions, which means that it is enough to assume that $u_0,\theta_0, d_0,\nabla d_0$ are H\"older continuous, choosing $\mu$ close to $1/p$ which is possible if $q$ is large enough.

Recalling  that the state manifold of (\ref{NLCF}) is given by
$$
\cSM=\{ (u,\theta, d)\in X_{\gamma}: \, \theta>0,\; |d|_2=1\},
$$
where $X_\gamma:=X_{\gamma,1}$, we see that by these results that $\cSM$ is locally positive invariant for the semi-flow, the total energy ${\sf E}$ is preserved and the negative total entropy $-{\sf N}$ is a strict Lyapunov functional for the semi-flow on $\cSM$. This completes the proof of Theorem \ref{thm:1}.

\rightline{$\Box$}

\vspace{.1cm}\noindent
{\em Proof of Theorems \ref{thm:2} and \ref{thm:3}}. \\
The linearization of  the system \eqref{eq:el}, \eqref{eq:thermo}, \eqref{eq:const} at an equilibrium $v_*=(0,\theta_*,d_*)$ is given by the operator
$A_*=A(v_*)$ defined in the base space $X_0$ with domain $D(A_*)=X_1$. This  operator has maximal $L_p$-regularity. Moreover,
$A_*$ is the negative generator of a compact analytic $C_0$-semigroup having compact resolvent, due to the compact embedding of $X_1=D(A_*)$ into $X_0$.
Hence, its spectrum consists only of countably many eigenvalues of finite multiplicity.

\begin{lemma}
Let  $z \ne 0$ be an eigenvalue of $A_*$. Then ${\rm Re}\,z < 0$.
\end{lemma}

{\em Proof}.
Suppose that $z\in \C\backslash\{0\}$ is an eigenvalue of $A_*$ with ${\rm Re}\, z\geq 0$. Then
\begin{align*}
\cL_\pi(z,\nabla)v_\pi &= 0\quad \mbox{ in } \Omega,\\
u=\partial_\nu \theta=\partial_\nu d &=0\quad \mbox{ on } \partial\Omega.
\end{align*}
where $(v_\pi)=(u,\theta,d)$ as above. Multiplying the equation for $d$ with $\bar{z}$ and taking the inner product of this equation with $v_\pi$ in $L_2(\Omega)$ yields by  integration by parts
the estimate
\begin{align*}
0={\rm Re}(\cL_\pi(z,\nabla)v_\pi|v_\pi)_{L_2} &\geq c\Big[ {\rm Re}\,z( |u|_{L_2}^2 + |\theta|_{L_2}^2 + |\nabla d|_{L_2}^2)\\
&\quad + |z|^2|d|_{L_2} + |\nabla u|_{L_2}^2 +|\nabla \theta|_{L_2}^2\Big].
\end{align*}
This implies $u=\theta=d=0$. Hence $A_*$ does not have eigenvalues in the $L_2$-setting with nonnegative real parts, except for $z=0$.
Due to elliptic regularity, eigenvalues are independent of $p$, and so the assertion follows for $A_*$ defined in $X_0$.

\rightline{$\Box$}

The above lemma states that all eigenvalues of $A_*$ expect for $0$ are stable. In addition, the eigenvalue $0$ is semi-simple. Its
eigenspace is given by
$$
{\sf N}(A_*)=\{(0,\vartheta, {\sf d}):\, \vartheta\in\RR, {\sf d}\in\RR^n\},
$$
and hence coincides with the set of constant equilibria $\bar{\cE}$ determined in Theorem \ref{prop:equi} when ignoring the constraint $|d|_2=1$ and
conservation of energy. Therefore each such equilibrium is normally stable.
Hence, the assertion of Theorem \ref{thm:2} follows by means of the generalized principle of linearized stability, Theorem 3.3.

Finally, we note that the assertion of Theorem \ref{thm:3} follows from Theorem 3.4, as we have compact embeddings and $-{\sf N}$ serves as a strict Lyapunov functional.

\rightline{$\Box$}


\begin{thebibliography}{10}


\bibitem{ASS04}
R. Atkin, T. Sluchin, I.W. Stewart,
\emph{Reflections on the life and work of Frank Matthews Leslie},
J. Non-Newtonian Fluid Mech. \textbf{119}, (2004), 7-23.

\bibitem{BP07}
D. Bothe, J. Pr\"uss,
\emph{$L^p$-theory for a class of non-Newtonian fluids},
SIAM J. Math. Anal. \textbf{39}, (2007), 379-421.


\bibitem{Cha92} S. Chandrasekhar,
Liquid Crystals, Cambridge University Press, (1992).

\bibitem{CS01}
D. Coutand, S. Shkoller,
\emph{Well-posedness of the full {E}ricksen-{L}eslie model of nematic liquid crystals},
{C. R. Acad. Sci. Paris S\'er. I Math.} \textbf{333}, (2001), 919-924.




\bibitem{DGP95}
P.G. DeGennes, J. Prost,
\emph{The Physics of Liquid Crystals}, Oxford University Press, (1995).


\bibitem{DHP03}
R. Denk, M. Hieber, J. Pr\"uss,
\emph{$\mathcal{R}$-boundedness, Fourier multipliers and problems of elliptic and parabolic type}.
Mem. Amer. Math. Soc., Vol. 166,  (2003).


\bibitem{Eri62}
J.~L. Ericksen, \emph{Hydrostatic theory of liquid crystals}, Arch. Rational
  Mech. Anal. \textbf{9} (1962), 371--378. 


\bibitem{FRS11}
E. Feireisl, E. Rocca, G. Schimperna,
\emph{On a non-isothermal model for nematic liquid crystals}, Nonlinearity, \textbf{24}, (2011), 243-257.

\bibitem{FFRS12}
E. Feireisl, M. Fr{\'e}mond, E. Rocca, G. Schimperna ,
\emph{A new approach to non-isothermal models for nematic liquid crystals},
Arch. Ration. Mech. Anal., \textbf{205}, (2012), 651-672.
	
\bibitem{FSRZ12}
E. Feireisl, E. Rocca, G. Schimperna, A. Zarnescu,
\emph{Evolution of non-isothermal Landau-de Gemmes nematic liquid crystal flows with singular potential},
arXiv:1207.1643

\bibitem{FSRZ13}
E. Feireisl, E. Rocca, G. Schimperna, A. Zarnescu,
\emph{Nonisothermal nematic liquid crystal flows with Ball-Majumdar free energy},
arXiv:1310.8474

\bibitem{Fra58}
F.C. Frank,
\emph{On the theory of liquid crystals}, Discussions Faraday Soc. \textbf{25}, (1958), 19-28.



\bibitem{HNPS14}
M. Hieber, M. Nesensohn, J. Pr\"uss, K. Schade,
\emph{Dynamics of nematic lquid crystals: the quasilinear approach},
Ann. Inst. H. Poincar\'e Anal. Non Lin\'eaire {\bf 33}, (2016),  397--408.

\bibitem{HP15a}
M. Hieber, J. Pr\"uss,
\emph{Thermodynamic Consistent Modeling and Analysis of Nematic Liquid Crystal Flows},
In: Springer Proc. Math. \& Statistics, 2016, to appear.

\bibitem{HP15b}
M. Hieber, J. Pr\"uss,
\emph{Modeling and Analysis of the Ericksen-Leslie equations for nematic liquid crystal flows},
In: Handbook of Mathematical Analysis in Mechanics of Viscous Fluids, (Y.Giga, A. Novotny (Eds.)), Springer, submitted.


\bibitem{HLX13}
M. Hong, J. Li, Z. Xin,
\emph{Blow-up criteria of strong solutions to the Ericksen-Leslie system in $\R^3$}.
Comm. Partial Differential Equations, \textbf{39}, (2014), 1284-1328.


\bibitem{HW10}
X. Hu, D. Wang,
\emph{Global solutions to the three-dimensional incompressible flow of liquid crystals}.
Comm. Math. Phys. \textbf{296}, (2010), 861-880.




\bibitem{HLLW15}
T. Huang, F. Lin, C. Liu, C. Wang,
\emph{Finite time singularities of the nematic liquid crystal flow in dimension three}. 
Arch. Rational Mech. Anal. \textbf{221}, (2016), 1223-1254.



\bibitem{HLW13}
J. Huang, F. Lin, C. Wang,
\emph{Regularity and existence of global solutions to the Ericksen-Leslie system in $\R^2$},
Comm. Math. Phys. \textbf{331}, (2014), 805-850. 

\bibitem{KPW10}
M.~K\"ohne, J.~Pr\"uss, M.~Wilke, \emph{On quasilinear parabolic evolution
  equations in weighted $L_p$-spaces}, J.~Evol.~Equ. \textbf{10} (2010), 443-463.

\bibitem{LPW14}
J. LeCrone, J. Pr\"uss, M. Wilke,
\emph{On quasilinear parabolic evolution equations in weighted $L^p$-spaces II.}
J. Evol. Equ. \textbf{14}, (2014), 509-533.


\bibitem{Les68}
F.~M. Leslie, \emph{Some constitutive equations for liquid crystals}, Arch.
  Rational Mech. Anal. \textbf{28} (1968), 265--283. 



\bibitem{Lin89}
F. Lin,
\emph{Nonlinear theory of defects in nematic liquid crystals: phase transition and flow phenomena}, Comm. Pure Appl. Math. \textbf{42}
  (1989), 789--814. 

\bibitem{Lin91}
F. Lin,
\emph{On nematic liquid crystals with variable degree of freedom}, Comm. Pure Appl. Math. \textbf{44}
  (1991), 453-468. 

\bibitem{LL95}
F. Lin, Ch. Liu, \emph{Nonparabolic dissipative systems modeling the flow of liquid crystals},
Comm. Pure Appl. Math. \textbf{48} (1995), 501--537. 

\bibitem{LL00}
F. Lin, Ch. Liu, \emph{Existence of solutions for the Ericksen-Leslie system},
Arch. Ration. Mech. Anal. \textbf{154} (2000), 135-156.

\bibitem{LLW10}
F. Lin, J. Lin, C. Wang,
\emph{Liquid crystal flows in two dimensions}.
Arch. Ration. Mech. Anal.  \textbf{197}, (2010), 297-336.

\bibitem{LW14b}
F. Lin, C. Wang,
\emph{Recent developments of analysis for hydrodynamic flow of nematic liquid crystals}.
Philos. Trans. R. Soc. Lon. Ser. A, Math. Phys. Eng. Sci \textbf{372}, (2014), 20130361.






\bibitem{MLG14}
W. Ma, H. Gong, J. Li,
\emph{Global strong solutions to incompressible Ericksen-Leslie system in $\R^3$}.
Nonlinear Anal. \textbf{109}, (2014), 230-235.


\bibitem{Ose33}
{C. W. Oseen},
\emph{The theory of liquid crystals}, Trans. Faraday Soc. \textbf{29}, (1933), 883-899.

\bibitem{Par70}
O. Parodi,
\emph{Stress tensor for a nematic liquid crystal}, J. Physique \textbf{31}, (1970), 581-584.

\bibitem{Pru03}
J. Pr{\"u}ss,
\emph{Maximal regularity for evolution equations in {$L_p$}-spaces},
Conf. Semin. Mat. Univ. Bari (2002), no.~285, (2003), 1-39.

\bibitem{PS04}
J.~Pr{\"u}ss, G.~Simonett, \emph{Maximal regularity for evolution equations
  in weighted $L_p$-spaces}, Arch.~Math.(Basel) \textbf{82} (2004),
  415--431.

\bibitem{PS16}
J.~Pr{\"u}ss, G.~Simonett, \emph{Moving Interfaces and Quasilinear Parabolic Evolution Equations},
Monographs of Mathematics {\bf 105}, Birkh\"auser, 2016.

\bibitem{PSZ09}
J.~Pr{\"u}ss, G.~Simonett, R.~Zacher, \emph{On convergence of solutions to
  equilibria for quasilinear parabolic problems}, J.~Diff.~Eqns. \textbf{246}
  (2009), 3902--3931.

\bibitem{Ste04}
I. W. Stewart, \emph{The Static and Dynamic Continuum Theory of Liquid Crystals}.
The Liquid Crystal Book Series, Taylor and Francis, 2004.


\bibitem{SL09}
H. Sun, Ch. Liu,
\emph{On energetic variational approaches in modeling the nematic liquid crystal flows},
Disc. Cont. Dyn. Syst. \textbf{23} (2009), 455-475.

\bibitem{Tri92}
H.~Triebel, \emph{Theory of Function Spaces II},
Birkh\"auser, Basel 1992

\bibitem{Vir94}
E. G. Virga,
\emph{Variational Theories for Liquid Crystals},
Chapman-Hall, London, 1994.


\bibitem{Wan11}
C. Wang,
\emph{Well-posedness for the heat flow of harmonic maps and the liquid crystal flow with rough initial data},
Arch. Ration. Mech. Anal. \textbf{200} (2011), 1-19.

\bibitem{WZZ13}
W. Wang, P. Zhang, Z. Zhang,
\emph{Well-posedness of the Ericksen-Leslie system},
Arch. Ration. Mech. Anal. \textbf{210} (2013), 837-855.

\bibitem{WXL13}
H. Wu, X. Xu, Ch. Liu,
\emph{On the general Ericksen-Leslie Ssystem: Parodi's relation, well-posedness and stability},
Arch. Ration. Mech. Anal. \textbf{208} (2013), 59-107.


\end{thebibliography}
\end{document}